\documentclass[11pt,a4paper,twoside,leqno]{amsart}
\usepackage[dvipsnames]{xcolor}
\usepackage{eucal}
\definecolor{britishracinggreen}{rgb}{0.0, 0.26, 0.15}
\definecolor{cobalt}{rgb}{0.0, 0.28, 0.67}
\definecolor{caribbeangreen}{rgb}{0.0, 0.8, 0.6}
\definecolor{brightpink}{rgb}{1.0, 0.0, 0.5}

\usepackage{tikz-cd}
\usepackage{amsmath}
\usepackage{amscd}
\usepackage{amsbsy}
\usepackage{amssymb}

\usepackage[english]{babel}
\usepackage{mathrsfs}
\usepackage{amsmath,caption,comment,mathtools,stmaryrd,enumitem}
\usepackage{multirow,booktabs,relsize}
\usepackage{fullpage}
\usepackage[colorlinks,bookmarks]{hyperref}
\usepackage{amsthm}

\usepackage{colonequals}

\theoremstyle{plain}
\newtheorem{theorem}{Theorem}[section]

\newtheorem{corollary}[theorem]{Corollary}
\newtheorem{proposition}[theorem]{Proposition}

\newtheorem{question}[theorem]{Question}

\theoremstyle{remark}

\newtheorem{remark}[theorem]{Remark}

\hypersetup{
    colorlinks,
    citecolor=britishracinggreen,
    filecolor=black,
    linkcolor=cobalt,
    urlcolor=black
}

\DeclareSymbolFont{usualmathcal}{OMS}{cmsy}{m}{n}
\DeclareMathAlphabet\BCal{OMS}{cmsy}{b}{n}

\DeclareSymbolFont{cmarrows}{OMS}{cmsy}{m}{n}
\SetSymbolFont{cmarrows}{bold}{OMS}{cmsy}{b}{n}
\DeclareMathSymbol{\cmminus}{\mathbin}{cmarrows}{"00}

\DeclareMathSymbol{\leftrightarrow}{\mathrel}{cmarrows}{"24}
\DeclareMathSymbol{\leftarrow}{\mathrel}{cmarrows}{"20}

\DeclareMathSymbol{\rightarrow}{\mathrel}{cmarrows}{"21}
\let\to=\rightarrow
\DeclareMathSymbol{\mapstochar}{\mathrel}{cmarrows}{"37}

\numberwithin{equation}{section}
\title{A note on specializing semi-orthogonal decompositions}
\begin{document}
\vspace{-3em}

\author[W.Wu]{Weimufei Wu}
\address{Peking University}
\email{2100010645@stu.pku.edu.cn}

\begin{abstract}
    We prove that families of semi-orthogonal decompositions do not satisfy the existence part of the valuative criterion for properness, giving a negative answer to Question~8.8 posed by Belmans, Okawa, and Ricolfi in \cite{Belmans2020ModuliSO}.
\end{abstract}
\maketitle
\section{Introduction}
Throughout, we work over the complex numbers $\mathbb C$. Let X be a smooth variety, and denote the bounded derived category of coherent sheaves on $X$ by $\mathrm {D^b}\mathrm{Coh}(X)$. A \textit{semi-orthogonal decomposition} of $\mathrm {D^b}\mathrm{Coh}(X)$ is a collection $(\mathcal{A}_1,\ldots,\mathcal{A}_n)$ of full triangulated subcategories of~$\mathrm {D^b}\mathrm{Coh}(X)$ such that
$\mathrm{Hom}_{\mathrm {D^b}\mathrm{Coh}(X)}(A_i, A_j)=0$ for all $A_i \in \mathcal{A}_i, A_j \in \mathcal{A}_j,j<i$, and that the smallest triangulated subcategory of $\mathrm {D^b}\mathrm{Coh}(X)$ containing $\mathcal{A}_1,\ldots,\mathcal{A}_n$ is $\mathrm {D^b}\mathrm{Coh}(X)$. We write $$\mathrm {D^b}\mathrm{Coh}(X)=\langle\mathcal{A}_1,\ldots,\mathcal{A}_n\rangle.$$

Given a smooth morphism of smooth schemes $X\to B$, a semi-orthogonal decomposition
    $$\mathrm {D^b}\mathrm{Coh}(X)=\langle\mathcal{A}_1,\ldots,\mathcal{A}_n\rangle$$
    is said to be $B$-linear if $\mathcal A_i\otimes f^*m\subset \mathcal A_i$ for any $m\in \mathrm {D^b}\mathrm{Coh}(B),$ $i=1,\ldots,n.$

In this short note, we prove the following result about specializations of semi-orthogonal decompositions.
\begin{theorem}\label{1.1}
  There exists a smooth projective family of rational surfaces $X\to \mathrm{Spec}( R)$ ($R$ is a discrete valuation ring, and $K$ is its field of fractions) with a $K$-linear semi-orthogonal decomposition of $\mathrm {D^b}\mathrm{Coh}(X_K)$ denoted by $\langle \mathcal{A}_K, \mathcal{B}_K\rangle$, which satisfies the following condition:

  \begin{itemize}
      \item there exists no $R$-linear semi-orthogonal decomposition of $\mathrm {D^b}\mathrm{Coh}(X_R)$ whose base change along $\mathrm{Spec} (K)\!\to\! \mathrm{Spec} (R)$ is $\langle\mathcal{A}_K, \mathcal{B}_K\rangle$. The base change for semi-orthogonal decompositions is in the sense of the base change for semi-orthogonal decompositions in \cite{Kuznetsov_2011}.
  \end{itemize}
\end{theorem}
This yields a counter-example to Question \ref{valuative} \cite[Question 8.8]{Belmans2020ModuliSO}, see Section 3.

\subsection{Acknowledgements}
This work was completed as part of the author's undergraduate research program at Peking University. The author is grateful to Pieter Belmans, Chunyi Li, and Shinnosuke Okawa for their helpful comments. The author would like to thank his supervisor Qizheng Yin for carefully reading this note and providing detailed guidance on the writing.

\section{Preliminaries}

We begin with the following theorem, which is a summary of how semi-orthogonal decompositions behave in families and under base changes in terms of moduli spaces.

\begin{theorem}[{\cite[Theorem A]{Belmans2020ModuliSO}}]\label{sod_f}
     Let $f\colon X\to B$ be a smooth projective morphism between smooth varieties. There exists an algebraic space $\mathrm{SOD}_f\to B$ which is \'etale over B and a functorial bijection
    $$\{B\textup{-morphisms }V\to \mathrm{SOD}_f \}\cong \{V\textup{-linear SODs } \mathrm {D^b}\mathrm{Coh}(X_V)=\langle\mathcal{A},\mathcal{B}\rangle\}$$
    for every quasi-compact and semi-separated $B$-scheme $V$.

\end{theorem}

In particular, this theorem implies:
\begin{enumerate}
    \item we have a well-defined pullback of semi-orthogonal decompositions linear over the base: When $v\colon V\to U$ is a morphism of $B$-schemes, we can pull $\mathrm {D^b}\mathrm{Coh}(X_U)=\langle\mathcal{A}_U,\mathcal{B}_U\rangle$ along~$v$ to get a $V$-linear semi-orthogonal decomposition $\mathrm {D^b}\mathrm{Coh}(X_V)=\langle\mathcal{A}_V,\mathcal{B}_V\rangle$. We will always write base changes of a semi-orthogonal decomposition in this way.
    \item we can extend a semi-orthogonal decomposition to an \'etale neighborhood (see \label{corB} \cite[Corollary B]{Belmans2020ModuliSO}). If there exists a semi-orthogonal decomposition $\mathrm {D^b}\mathrm{Coh}(X_p)=\langle\mathcal{A},\mathcal{B}\rangle$ where $X_p$ is a fiber over a closed point $p\colon\mathrm{Spec}(\mathbb C)\to B$ , then there exists an \'etale neighborhood $U$ of $p$ ($p$ factors through $U$)  such that:

\end{enumerate}

\begin{itemize}
    \item there exists a unique semi-orthogonal decomposition $\mathrm {D^b}\mathrm{Coh}(X_U)=\langle\mathcal{A}_U,\mathcal{B}_U\rangle$ whose pullback to $p$ is $\mathrm {D^b}\mathrm{Coh}(X_p)=\langle\mathcal{A},\mathcal{B}\rangle.$
\end{itemize}

We will also need a weak version of the theorems on $K$-theories of linear categories proved in \cite{Blanc_2016}, \cite{moulinos2019derived}, and \cite{Perry_2022}. We refer the reader to Section 2 in \cite{Perry_2022} for the formalism of categories linear over a base scheme. For our purposes, we only need to know that: for a morphism of smooth schemes $f\colon X \to B$, we have an induced structure of a $B$-linear category on  $\mathrm{D^b}\mathrm{Coh}(X)$; a $B$-linear semi-orthogonal decomposition of $\mathrm{D^b}\mathrm{Coh}(X)$ is a semi-orthogonal decomposition of $\mathrm{D^b}\mathrm{Coh}(X)$ as a $B$-linear category, and every component in this semi-orthogonal decomposition is also a $B$-linear category.

\begin{theorem}[see {\cite[Section 5]{Perry_2022}}]\label{Ktheory}
    There exists a relative topological $K$-theory functor from
    $B$-linear categories to abelian groups $$K_0^\mathrm{top}(-/B): {\mathrm {Cat}_B \to \textup{Abelian groups}}$$ such that:
    \begin{enumerate}
        \item for $X$ a variety over
 $B=\mathrm{Spec} ( \mathbb C)$, we have $K^\mathrm{top}_0(\mathrm {D^b}\mathrm{Coh}(X)/B)=K^\mathrm{top}_0(X^\mathrm{an})$,
where the right hand side is the complex $K$-theory group of the analytification of $X$.
\item for a $\mathbb C$-linear category $\mathcal{C} $ there exists a functorial map $K_0(\mathcal{C})\to K_0^\mathrm{top}(\mathcal{C}/\mathbb C)$, where $K_0$ is the Grothendieck group of a $\mathbb C$-triangulated category;
 \item if $\mathcal C = \langle \mathcal C_1, \ldots, \mathcal C_m\rangle$ is a semi-orthogonal decomposition of $B$-linear categories, denote the corresponding inclusions by $\iota_i:\mathcal{C}_i\to \mathcal C$ and the projections by $\iota_i^*:\mathcal{C}\to \mathcal C_i$. These functors induce a splitting$$K^\mathrm{top}_0(\mathcal C/B)\cong \bigoplus_{i=1}^m K^\mathrm{top}_0(\mathcal C_i/B),$$
and when $B=\mathrm{Spec} ( \mathbb C)$, these functors also yield $K_0(\mathcal C)\cong \bigoplus_{i=1}^m K_0(\mathcal C_i)$;
 \item  if $f\colon X \to  B$ is a smooth proper morphism of varieties, $\mathcal{C}\subset \mathrm {D^b}\mathrm{Coh}(X)$ is a $B$-linear admissible subcategory, then $K^\mathrm{top}_0(\mathcal C/B)$ is a local system of finitely generated abelian groups on $B$ whose fiber over any point $p \in B(\mathbb C)$ is $K^\mathrm{top}_0(\mathcal C _p)$.
    \end{enumerate}
\end{theorem}
\section{Proof of the main result}
To prepare for our construction, we first prove a result about the deformation of semi-orthogonal decompositions.
\begin{proposition}\label{phantom spread}
   Let~$X\to B$ be a smooth projective morphism between smooth connected varieties, and $\mathrm {D^b}\mathrm{Coh}(X)=\langle\mathcal{A},\mathcal{B}\rangle$ a $B$-linear semi-orthogonal decomposition.
For a fiber $X_b$ over a point $b$ with a full exceptional collection $\mathrm {D^b}\mathrm{Coh}(X_b) =\langle E_1,\ldots,E_n\rangle$, we have $K_0(\mathcal{A}_{c})=K_0(\mathcal{A}_b)$ for any closed point $c$ in an \'etale open neighborhood of $b$.
\end{proposition}

\begin{proof}
    Since the local deformation of an exceptional collection will also be an exceptional collection, (see \cite[Remark 7.6]{Belmans2020ModuliSO}), we can shrink $B$ to ensure that $\mathrm {D^b}\mathrm{Coh}(X)$ has a relative full exceptional collection, still denoted by $\langle E_1,\ldots,E_{n}\rangle$.
Now we apply (2) and (3) of Theorem~\ref{Ktheory} to the two semi-orthogonal decompositions of $\mathrm {D^b}\mathrm{Coh}(X_b)$: $\mathrm {D^b}\mathrm{Coh}(X_b)=\langle E_1, \ldots,E_{n}\rangle$ and $\mathrm {D^b}\mathrm{Coh}(X_b)=\langle\mathcal{A}_b,\mathcal{B}_b\rangle$. Then we get two commutative diagrams of groups. 
$$
\begin{tikzcd}
\bigoplus K_0(E_i) \arrow[r, "\cong"] \arrow[d, "\cong"] & \bigoplus K_0^{\mathrm{top}}(E_i) \arrow[d, "\cong"] \\
K_0(X_b) \arrow[r]                                       & K_0^{\mathrm{top}}(X_b)                             
\end{tikzcd}    \hspace{30pt}                             
\begin{tikzcd}
K_0(\mathcal{A}_b) \arrow[r] \arrow[d, hook] & K_0^{\mathrm{top}}(\mathcal{A}_b) \arrow[d, hook] \\
K_0(X_b) \arrow[r, "\cong"]                  & K_0^{\mathrm{top}}(X_b)                          
\end{tikzcd}$$
Note the comparison morphism $K_0(\langle E_i\rangle_b)\to K^\mathrm{top}_0(\langle E_i\rangle_b)$ is an isomorphism for $i=1,\ldots,n$, so we have an induced isomorphism $K_0(X_b)\to K^\mathrm{top}_0(X_b)$ by the left diagram. Since the left and right arrows in the right diagram are both injective by (3) of Theorem \ref{Ktheory}, we obtain that the above arrow in the right diagram $K_0(\mathcal{A}_b)\to K_0^\mathrm{top}(\mathcal{A}_b)$ is an isomorphism.

For any closed point $c$, we have $$K_0(\mathcal{A}_c)=K^\mathrm{top}_0(\mathcal{A}_c)=K^\mathrm{top}_0(\mathcal{A}_b)=K_0(\mathcal{A}_b).$$ The first identity follows the same line as the argument about $\mathcal{A}_b$ and the second identity follows from the fact that $K^\mathrm{top}_0(\mathcal{A}/B)$ is a local system by (4) of Theorem \ref{Ktheory}.
\end{proof}

For the case we are interested in, we recall the definition of a phantom category. A nonzero admissible subcategory $\mathcal A \subset  \mathrm {D^b}\mathrm{Coh}(X)$ is said to be a \textit{phantom} if its Grothendieck group $K_0(\mathcal A)$ is trivial; see \cite{phantom ihes}.  

\begin{corollary}
    If $\mathcal{A}_b$ is a phantom, then $K_0(\mathcal{A}_c)=0$.
\end{corollary}
We use the phantom constructed by Krah in \cite{Krah} to  prove Theorem \ref{1.1}, which gives a negative answer to the following question posed in \cite{Belmans2020ModuliSO}.

\begin{question}[{\cite[Question 8.8]{Belmans2020ModuliSO}}]\label{valuative}
   We keep the assumptions of Theorem \ref{sod_f}. Let $R$ be a discrete valuation ring with field of fractions $K$. For every commutative diagram as follows, does there exist a dotted arrow that makes the diagram commutative?
   $$\begin{tikzcd}
\mathrm{Spec} (K) \arrow[r] \arrow[d]         & \mathrm{SOD}_f \arrow[d] \\
\mathrm{Spec}(R) \arrow[r] \arrow[ru, dotted] & B              
\end{tikzcd}$$
\end{question}

\begin{theorem}[\cite{Krah}, Theorem 1.1]\label{krah}
    Let $X$ be the blow-up of $\mathbb P^2_\mathbb C$ at ten closed points in general position. Then there exists a phantom in $\mathrm {D^b}\mathrm{Coh}(X)$, which is the \textit{right-orthogonal complement} of an exceptional collection of length $13$.
\end{theorem}
To construct our example, we choose a smooth projective morphism $f\colon\mathcal X\to B\cong \mathbb P^1$ such that: fix a point $0\in B$, the fiber $\mathcal{X}_0 \colonequals f^{-1}(0)$ is a blow-up of $\mathbb P^2$ at $10$ points whose positions are to be determined; there exists an open subset $B^o\subset B$ such that for any point $u\in$$B^o$, the fiber $\mathcal X_u$ is a blow-up of $\mathbb P^2$ in $10$ general points. 

We claim that after suitable modification, $f$ provides a counter-example to Question \ref{valuative}.
We fix a point $p\in B^o$. By Theorem \ref{krah}, $\mathcal{X}_p$ has a semi-orthogonal decomposition $\langle\mathcal{B},\mathcal{P}\rangle$ with $\mathcal{P}$ being a phantom.
 By (\ref{corB}) under Theorem \ref{sod_f}, we can replace $B$ with a finite cover such that:
 \begin{enumerate}
     \item there exists an open subscheme $B^o\subset B $ such that there exists a $B^o$-linear semi-orthogonal decomposition $\langle\mathcal{B}',\mathcal{P}'\rangle$ of the category $\mathrm{D^b Coh}(\mathcal X\times_B B^o)$;
     \item the base change of $\mathcal{P}'$ to $\mathcal X_p$ coincides with $\mathcal{P}$.
 \end{enumerate}
 
 To satisfy these conditions, we take an \'etale neighborhood $B_1$ of $p$ and then take a compactification $\bar B_1$ such that $0$ has an inverse image in $\bar B_1 \to B$; then we replace $B$ with $\bar B_1$, and replace $0$ with one of its inverse images.
Let $\langle\mathcal{B}_K,\mathcal{P}_K\rangle$ denote the pullback of $\langle\mathcal{B}',\mathcal{P}'\rangle$ along $\mathrm{Spec} (K)\to B^o$. Denote the local ring of $B$ at $0$ by $(R,k)$, where $k$ is the residue field.

\begin{proof} [Proof of Theorem \ref{1.1}]
\renewcommand{\qedsymbol}{}
    Our goal is to prove that we cannot specialize $\langle\mathcal{B}_K,\mathcal{P}_K\rangle$, i.e. there exists no $R$-linear semi-orthogonal decomposition whose base change to $\mathrm{Spec}(K)$ is $\langle\mathcal{B}_K,\mathcal{P}_K\rangle$. We will prove this by contradiction. Suppose there exists such a decomposition $\mathrm {D^b}\mathrm{Coh}(\mathcal X_R)=\langle\mathcal{B}_R,\mathcal{P}_R\rangle$. Firstly, we pull $\langle\mathcal{B}_R,\mathcal{P}_R\rangle$ along $\mathrm{Spec} (k)\to \mathrm{Spec}( R)$ to $\langle\mathcal{B}_k,\mathcal{P}_k\rangle$.
\end{proof}
\begin{proposition}
    $ \mathcal{P}_k$ is a phantom. 
\end{proposition}
\begin{proof}
We can extend $\mathcal{P}_k$ by (\ref{corB}) under Theorem \ref{sod_f} to get a connected elementary \'etale open neighborhood $C\to \mathrm{SOD}_f$. Since $C$ and $\mathrm{SOD}_f$ are both \'etale over $B$, we know that $C\to \mathrm{SOD}_f$ is also \'etale. After an \'etale base change, we can assume that $\mathrm{Spec} (R)$ factors through $C$. Hence the fiber product $C\times_{\mathrm{SOD}_f}B^o$ is nonempty.

We now apply Proposition \ref{phantom spread} to $C $ and $B^o$ respectively to obtain $ K_0(\mathcal{P}_k)=0$. So if $\mathcal{P}_k$ is not a phantom, it can only be zero. But since we know that the local deformation is unique, we can choose $C$ so that $\mathcal{P}_C$ is trivial. However, passing to $C\times_ {\mathrm{SOD}_f}B^o$ will induce $\mathcal{P}_b=0$ for a general point $b\in B^o$. But this cannot happen because
\begin{itemize}
    \item there exists a $B^o$-linear semi-orthogonal decomposition $\langle\mathcal{B}^o,\mathcal{P}^o\rangle$ which coincides with $\langle\mathcal{B},\mathcal{P}\rangle$ at $p$, and $\mathcal{P}^o_q\neq 0$ for any $q\in B^o$ (this follows from the proof of Krah's example);
    \item the local deformation of $\mathcal{B}$ is unique, so for a general point $b\in B^o$, $\mathcal{P}_b=\mathcal{P}^o_b.$
\end{itemize}
So the corresponding $\mathcal{P}_b$ is not trivial, which is a contradiction.
\end{proof}
We can choose $\mathcal X_0$ to be the variety described in the following theorem. This concludes our proof of Theorem~\ref{1.1}.
\qed

\begin{theorem}[{\cite[Theorem 1.1]{Borisov2024NonexistenceOP}}]\label{lev}
     Let $X$ be the blow-up of $\mathbb P^2_\mathbb C$ at a finite number of points in a very general position on a smooth cubic curve. Then X has no phantom.
\end{theorem}
\begin{remark}
    If $X$ is a blow-up of $\mathbb P^2_\mathbb C$ in at most $9$ points and admits a phantom, denoted by $\mathcal{P} $, then this phantom extends to an \'etale neighborhood $U$. Now $\mathcal{P}$ is smooth proper over $U$, but we have $\mathcal P_\eta=0$, where $\eta$ is the generic point of $U$. This is because we can deform $X$ to a del Pezzo surface and this phantom deforms either to $0$ or to a phantom on a del Pezzo surface, but it must be zero by \cite[Theorem 6.35]{Pirozhkov2020AdmissibleSO}. So a consequence of \cite[Conjecture 8.29]{Belmans2020ModuliSO} is that~$X$ cannot admit a phantom.
    
\end{remark}

   \begin{remark}
 There are several natural questions related to specializing semi-orthogonal decompositions. The second question follows the idea in Kulikov models of Type I.
  
\begin{enumerate}
        \item  Is there an explicit algebraic obstruction to specializing semi-orthogonal decompositions?

        \item  Denote the semi-orthogonal decomposition corresponding to $\mathrm{Spec}(K)\to\mathrm{SOD}_f$ in Question \ref{valuative} by $\langle\mathcal{A}_K,\mathcal{B}_K\rangle$. Is Question \ref{valuative} valid given that $\mathcal{A}_K$ is a K3 category? 
    \end{enumerate}
   \end{remark}

\end{document}